\documentclass[dvipsnames]{article}

\title{Non-distributive positive logic as a fragment of \\ 
       first-order logic over semilattices}
\author{Jim de Groot \bigskip \\ 
\footnotesize{The Australian National University} \vspace{-.4em} \\ 
\footnotesize{Canberra, Ngunnawal Country, Australia}}
\date{}

\usepackage{amsmath}
\usepackage{amssymb}
\usepackage{amsthm}
  \theoremstyle{definition}
  \swapnumbers
    \newtheorem{para}{}[section]
    \newtheorem{definition}[para]{Definition}
    \newtheorem{example}[para]{Example}
    
    \newtheorem{remark}[para]{Remark}
  \theoremstyle{theorem}
    \newtheorem{lemma}[para]{Lemma}
    \newtheorem{corollary}[para]{Corollary}
    \newtheorem{theorem}[para]{Theorem}
    \newtheorem{proposition}[para]{Proposition}

\usepackage{stmaryrd}
\usepackage{tikz-cd}
  \usetikzlibrary{arrows,arrows.meta}
  \tikzcdset{arrow style=tikz, diagrams={>=stealth'}}
  \usetikzlibrary{snakes}
  \tikzset{decoration={snake, amplitude=.2mm,segment length=2mm}}
\usepackage[a4paper,margin=2.5cm]{geometry}
\usepackage[colorlinks=true,
            linktocpage=true,
            citecolor=RubineRed,
            linkcolor=RoyalBlue]{hyperref}
\usepackage{enumitem}
\usepackage{float}
\usepackage[numbers]{natbib}
\usepackage{doi}

\usepackage{titlesec}
\titleformat{\section}
   {\normalfont\large\bfseries\centering}{\thesection}{1em}{}
     [\vspace{-.5em}\rule{12em}{.15mm}]
\titleformat{\subsection}
   {\normalfont\bfseries\centering}{}{0em}{}%
     [\vspace{-.5em}\rule{12em}{.15mm}]
\titlespacing*{\section} {0pt}{6ex plus 1ex minus .2ex}{2ex plus .2ex}

\usepackage{mathrsfs}
\usepackage[cal=txupr,scr=boondoxupr,scrscaled=1.050]{mathalfa}

\newcommand{\ms}[1]{\mathscr{#1}}
\newcommand{\cat}[1]{\mathsf{#1}}   
\newcommand{\mo}[1]{\mathfrak{#1}}  
\newcommand{\lan}[1]{\mathbf{#1}}   
\newcommand{\fun}[1]{\mathscr{#1}}
\newcommand{\fmo}[1]{\mathcal{#1}}

\renewcommand{\iff}{\quad\text{iff}\quad}
\renewcommand{\phi}{\varphi}
\newcommand{\und}[1]{\underline{#1}}

\DeclareMathOperator{\Prop}{Prop}
\renewcommand{\th}{\operatorname{th}}
\DeclareMathOperator{\st}{st}
\DeclareMathOperator{\MOC}{MOC}
\newcommand{\FOL}{\lan{FOL}}

\usepackage{stackengine}
\newcommand{\cvee}{\mathrel{\stackon[.5pt]{$\vee$}{\tiny{$c$}}}}

\newcommand{\simul}{\mathrel{\mkern1mu\underline{\mkern-1mu \to\mkern-2mu}\mkern2mu }}
\newcommand{\login}{\rightsquigarrow}

\newcommand{\ismeet}{\mathsf{ismeet}}

\DeclareMathOperator{\Gr}{Gr}
\DeclareMathOperator{\id}{id}

\newcommand{\FSL}{\cat{FSL}}
\newcommand{\fmeet}{\curlywedge}
\newcommand{\fleq}{\preccurlyeq}
\newcommand{\Lone}{L$_1$\kern-.2ex}
\newcommand{\llb}{\llbracket}
\newcommand{\rrb}{\rrbracket}

\begin{document}

\maketitle

\begin{abstract}
  \noindent
  We characterise non-distributive positive logic as the fragment of
  a single-sorted first-order language that is preserved by a new notion
  of simulation called a meet-simulation.
  Meet-simulations distinguish themselves from simulations because they
  relate pairs of states from one model to single states from another.
  En route to this result we use a more traditional notion of simulations
  and prove a Hennessy-Milner style theorem for it, using an analogue of
  modal saturation called meet-compactness.
\end{abstract}

\paragraph{Version of record.}
This article has been accepted for publication in the Journal of Logic and Computation Published by Oxford University Press, doi:\href{dx.doi.org/10.1093/logcom/exad003}{\color{RoyalBlue}{10.1093/logcom/exad003}}.


\section{Introduction}

  The celebrated Van Benthem characterisation theorem characterises normal
  modal logic as the bisimu\-la\-tion-invariant fragment of first-order logic
  \cite{Ben76}.
  More precisely, it states that a first-order formula with one free variable
  is equivalent to the standard translation of a modal formula if and only if
  it is invariant under bisimulations.
  Similar theorems have been proven for non-normal modal logics,
  such as monotone modal logic~\cite{Han03}, neighbourhood logic~\cite{HanKupPac09}, and
  instantial neighbourhood logic~\cite{Gro22-inl}.
  For the latter three we only consider equivalence over particular classes of
  first-order structures, namely those corresponding to the modal logic semantics.
  Apart from modal extensions of classical propositional logic,
  other propositional logics and their modal extensions can be
  viewed as fragments of first-order logic as well.
  Examples include (bi-)intuitionistic logic~\cite{Rod86,KurRij97,Bad16}
  and positive modal logic~\cite{KurRij97,CelJan12}.
  All such theorems are instances of (model-theoretic) preservation theorems,
  see for example~\cite[Section~5.2]{ChaKei73} or~\cite{Com83}.

  Sometimes
  ``invariant under bisimulations'' is replaced with
  ``preserved by simulations.''
  This allows one to exclude negations from the characterisation.
  Indeed, if a formula is invariant under bisimulations then so is its
  negation, but preservation by simulations does not imply preservation
  of its negation.

  The goal of the present paper is to characterise non-distributive positive 
  logic~\cite{LuiChi76,Urq78,ConEA21} as a fragment of first-order logic.
  It has recently been noted that non-distributive positive logic logic can
  be given frame semantics
  by means of meet-semilattices with a valuation, with filters serving as
  denotations of formulae~\cite{Dmi21,BezEA22}.
  Conjunctions are interpreted as usual, while 
  $\phi \vee \psi$ holds at a state $w$ if there are two states $v$ and $u$ satisfying
  $\phi$ and $\psi$, respectively, such that their meet lies below $w$.
  This non-standard interpretation of disjunctions prevents distributivity.
  
  Since meet-semilattices are in particular partially ordered sets, they can be viewed as
  interpreting structures for the first-order language with one binary predicate.
  Indeed, the partial order underlying a meet-semilattice then serves as the
  interpretation of the binary relation.
  Guided by this observation, we translate
  non-distributive positive logic into the first-order logic with one binary
  predicate and a unary predicate for each proposition letter.
  This translation is similar to the one in~\cite{BezEA22}.
  
  To establish this characterisation result, we have to define an
  appropriate analogue of simulations. The main challenge is to exclude
  classical (locally evaluated) disjunctions from being preserved by this
  notion. Ordinary simulations do not suffice because they preserve
  all locally evaluated monotone connectives.
  Instead, we introduce \emph{meet-simulations}.
  These differ from ordinary simulations because they relate pairs of
  states from one model to states of another.
  This modification, it turns out, allows us to prevent preservation
  of classical disjunctions. We prove the following characterisation theorem:
  \begin{center}
    A first-order formula $\alpha(x)$ is
    preserved by meet-simulations if and only if \\
    it is equivalent to a $\neg(x = x)$ or to
    the standard translation of a \\ non-distributive positive formula
    over the class of meet-semilattices.
  \end{center}
  Observe that we can replace $\neg(x = x)$ with any contradiction,
  such as $x \wedge \neg x$,
  so it is possible to work in a first-order language without equality.
  To eliminate reference to $\neg(x = x)$ entirely,
  we finally suggest an adaptation of meet-simulations to
  meet$_{\omega}$-simulations.
  The special feature of meet$_{\omega}$-simulations is that they relate
  finite subsets of states of one model to states of another.
  With this adaptation, we obtain:
  \begin{center}
    A first-order formula $\alpha(x)$ is
    preserved by meet$_{\omega}$-simulations if and only if \\
    it is equivalent to
    the standard translation of a non-distributive \\ positive formula
    over the class of meet-semilattices.
  \end{center}

\paragraph{Related work.}
  The idea of using (bi)simulations to relate pairs of states to pairs of states
  appeared before in~\cite{FleEA15}, where it is used to characterise
  fragments of the calculus of binary relations, and in~\cite{AbrDesFig17}
  in the context of paths on data trees.

\paragraph{Structure of the paper.}
  In Section~\ref{sec:2} we recall the language and semantics of
  non-distributive positive logic. In Section~\ref{sec:3} we summarise
  some basic first-order logic required for the results in this paper,
  and we define the standard translation for non-distributive
  positive logic.
  
  The contributions of this paper start in Section~\ref{sec:4},
  where we define simulations between models and show that states related
  by a simulation are logically inclusive (that is, the theory of the former
  is included in the theory of the latter).
  In Section~\ref{sec:5} we give a Hennessy-Milner style theorem for simulations.
  We first prove that the finite models form a Hennessy-Milner class.
  Taking stock of the proof, we define an analogue of modal saturation called
  meet-compactness. We then prove that $\omega$-saturated models are
  meet-compact, and that the meet-compact models form a Hennessy-Milner class.
  
  Subsequently, in Section~\ref{sec:towards-VB} we work towards a Van Benthem style
  characterisation and point out why simulations do not admit
  such a result. We then introduce meet-simulations in Section~\ref{sec:VB}
  and use the results from Section~\ref{sec:towards-VB} to prove the
  characterisation theorem announced above.
  Finally, in Section~\ref{sec:7b} we show how to adapt the notion of a meet-simulation
  to prevent preservation of classical contradictions.
  We conclude in Section~\ref{sec:8}.

\section{Non-distributive positive logic}\label{sec:2}

  Denote by $\lan{L}$ the language of positive logic, i.e.~the language
  generated by the grammar
  $$
    \phi ::= p \mid \top \mid \bot \mid \phi \wedge \phi \mid \phi \vee \phi,
  $$
  where $p$ ranges over some arbitrary but fixed set $\Prop$ of proposition letters.
  One can define a logic of consequence pairs of $\lan{L}$-formula
  whose algebraic semantics is given by lattices, see~\cite[Section~3.1]{BezEA22}.

  Formulae from $\lan{L}$ can be interpreted in
  meet-semilattices with a valuation~\cite{Dmi21}.
  The intuition behind this is that the collection of filters of a
  meet-semilattices is closed under (arbitrary) intersections.
  Therefore they form a complete lattice, but disjunctions are not
  given by unions. This gives rise to a non-standard interpretation of
  disjunctions which prevents distributivity.
  We call the resulting frames and models \Lone-frames,
  to distinguish them from the slightly different L-frames
  used in~\cite{BezEA22}, see Remark~\ref{rem:L-vs-L1}.
  
  In this paper, by a \emph{meet-semilattice} we mean a partially ordered
  set in which every finite subsets has a greatest lower bound, called its \emph{meet}.
  The meet of $w$ and $v$ is denoted by $w \fmeet v$, reserving
  the symbol $\wedge$ for conjunctions of formulae.
  The empty meet is the top element, denoted by $1$.
  If $(W, 1, \fmeet)$ is a meet-semilattice then we write $\fleq$ for the
  partial order given by $w \fleq v$ iff $w \fmeet v = w$.
  A \emph{filter} of a meet-semilattice $(W, 1, \fmeet)$ is a subset $F$ of $W$
  which is upward closed under $\fleq$ and closed under finite meets.
  Filters are nonempty because they contain the empty meet, $1$.

\begin{definition}
  An \emph{\Lone-frame} is a meet-semilattice $(W, 1, \fmeet)$.
  An \emph{\Lone-model} is an \Lone-frame $(W, 1, \fmeet)$ together with a valuation
  $V$ that assigns to each proposition letter $p \in \Prop$ a
  filter $V(p)$ of $(W, 1, \fmeet)$.
  The interpretation of $\phi \in \lan{L}$ at a state $w$ of
  an \Lone-model $\mo{M} = (W, 1, \fmeet, V)$ is defined recursively via
  \begin{align*}
    \mo{M}, w \Vdash p &\iff w \in V(p) \\
    \mo{M}, w \Vdash \top &\phantom{\iff}\text{always} \\
    \mo{M}, w \Vdash \bot &\iff w = 1 \\
    \mo{M}, w \Vdash \phi_1 \wedge \phi_2 &\iff \mo{M}, w \Vdash \phi_1
                     \text{ and } \mo{M}, w \Vdash \phi_2 \\
    \mo{M}, w \Vdash \phi_1 \vee \phi_2
      &\iff \exists u, v \in W \text{ s.t. } u \fmeet v \fleq w
            \text{ and } \mo{M}, u \Vdash \phi_1 \text{ and } \mo{M}, v \Vdash \phi_2
  \end{align*}
  We denote the \emph{truth set} of $\phi$ by
  $\llb \phi \rrb^{\mo{M}} := \{ w \in W \mid \mo{M}, w \Vdash \phi \}$.
  The \emph{theory} of a state $w \in W$ is defined by
  $\th_{\mo{M}}(w) := \{ \phi \in \lan{L} \mid \mo{M}, w \Vdash \phi \}$.
  If $w$ and $w'$ are states in $\mo{M}$ and $\mo{M}'$, respectively,
  then we write $\mo{M}, w \login \mo{M}', w'$ if
  $\th_{\mo{M}}(w) \subseteq \th_{\mo{M}'}(w')$.
\end{definition}

  It can be shown that the truth set of every formula in any \Lone-model $\mo{M}$
  is a filter in the underlying \Lone-frame.
  In particular, the truth set of a formula of the form $\phi_1 \vee \phi_2$
  is the smallest filter containing both $\llb \phi_1 \rrb^{\mo{M}}$
  and $\llb \phi_2 \rrb^{\mo{M}}$, and can also be given by
  $
      \{ w \in W \mid u \fmeet v \fleq w \text{ for some }
           u \in \llb \phi_1 \rrb^{\mo{M}} \text{ and }
           v \in \llb \phi_2 \rrb^{\mo{M}} \}.
  $

  The following notion of morphism is taken from~\cite[Definition~3.4.9]{Dmi21}.

\begin{definition}
  An \emph{\Lone-morphism} between \Lone-frames
  $(W, 1, \fmeet)$ and $(W', 1', \fmeet')$ is a function $f : W \to W'$
  that preserves all finite meets,
  and satisfies for all $w \in W$ and $u', v' \in W'$: 
  \begin{itemize}
    \item $f(w) = 1'$ iff $w = 1$;
    \item If $u' \wedge v' \fleq' f(w)$, then there exist $u, v \in W$ such that
          $u' \fleq' f(u)$ and $v' \fleq' f(v)$ and $u \fmeet w \fleq w$.
  \end{itemize}
  The second condition can be depicted as follows:
  $$
    
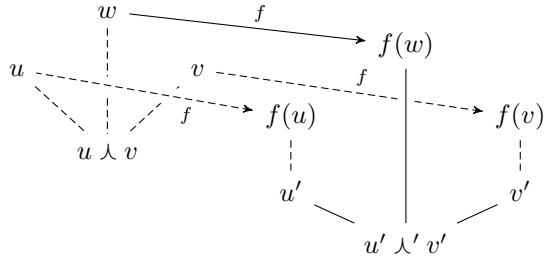
\begin{tikzcd}[scale=1.1]
        & [-1.3em]
          w \arrow[rrrd, "f"]
            \arrow[dddd, dashed, -]
        & [-1.3em]
        & [-1em]
        & [-1.5em]
        & [-1.5em] \\ [-2em]
        &
        &
        &
        & f(w)
        & \\ [-2.3em]
      u     \arrow[ddr, dashed, -]
        &
        & v \arrow[ddl, dashed, -]
            \arrow[rrrd, dashed, "f"]
        &&& \\ [-1.5em]
        &
        &
        & f(u) \arrow[dd, dashed, -]
               \arrow[lllu, dashed, <-, crossing over, "f" pos=.25]
        &
        & f(v) \arrow[dd, dashed, -] \\ [-2em]
        & u \fmeet v
        &&&& \\  [-1.7em]
        &
        &
        & u' \arrow[dr, -]
        &
        & v' \arrow[dl, -] \\ [-1.3em]
        &
        &
        &
        & u' \fmeet' v'  \arrow[uuuuu, crossing over, -]
        &
    \end{tikzcd}
  $$
  An \emph{\Lone-morphism} between \Lone-models $(W, 1, \fmeet, V)$ and
  $(W', 1', \fmeet', V')$ is an \Lone-morphism $f$ between the underlying
  frames such that $V(p) = f^{-1}(V'(p))$ for all $p \in \Prop$.
\end{definition}

\begin{remark}\label{rem:L-vs-L1}
  In~\cite{BezEA22} a variation of the semantics presented above is used where
  meet-semilattices are only assumed to have binary meets.
  As a consequence, they need not have a top element, and filters are allowed to be empty.
  An advantage of this approach is that the truth set of $\bot$ is empty,
  so that models contain no inconsistent state.
  The drawback is that it complicates interpretation of $\phi_1 \vee \phi_2$,
  cf.~\cite[Definition~3.6]{BezEA22}.
  This, in turn, affects the definition of truth-preserving morphisms between
  models.
  
  In this paper we choose to allow an inconsistent state $1$ satisfying $\bot$
  because it simplifies the definitions of simulations and meet-simulations.
\end{remark}

\section{First-order translation}\label{sec:3}

  We define the first-order language we work with and the standard translation of
  $\lan{L}$ into this language.
  Then we characterise the class of first-order structures corresponding to
  \Lone-models and we recall the definition of $\omega$-saturation.

\begin{definition}
  Let $\FOL$ be the single-sorted first-order language which has a
  unary predicate $P_p$ for every proposition letter $p \in \Prop$,
  and a binary predicate $R$.
  To avoid confusion with the interpretation of disjunctions from $\lan{L}$,
  we denote classical disjunctions like the one in $\FOL$ by $\cvee$.
\end{definition}

  Models for $\FOL$ are denoted by $\fmo{M}, \fmo{N}$ (as opposed to
  $\mo{M}, \mo{N}$ for \Lone-frames).
  Intuitively, the relation symbol of our first-order language accounts
  of the poset structure of \Lone-frames.
  If $x, y$ and $z$ are variables, then we can express that
  $x$ is the meet of $y$ and $z$ in the ordering induced by the
  interpretation of $R$ using a first-order sentence.
  To streamline notation we abbreviate this as follows:
  $$
    \ismeet(x, y, z)
      := (x R y) \wedge (x R z) \wedge \forall x'((x' R y \wedge x' R z) \to x' Rx).
  $$
  We are now ready to define the standard translation.

\begin{definition}\label{def:st}
  Let $x$ be a variable.
  Define the standard translation $\st_x : \lan{L} \to \FOL$ recursively
  via
  \begin{align*}
    \st_x(p) &= P_px \\
    \st_x(\top) &= (x = x) \\
    \st_x(\bot) &= \forall y(yRx) \\
    \st_x(\phi \wedge \psi) &= \st_x(\phi) \wedge \st_x(\psi) \\
    \st_x(\phi \vee \psi)
      &= \exists x' \exists y \exists z (\ismeet(x', y, z) \wedge x'Rx
         \wedge \st_y(\phi) \wedge \st_z(\psi))
  \end{align*}
\end{definition}

  While one often sees the standard translation of $\top$ defined as $x = x$,
  it is not strictly necessary: any tautology
  with free variable $x$ suffices.
  Thus, we can also work in a first-order language without equality.

  Every \Lone-model $\mo{M} = (W, 1, \fmeet, V)$ gives rise to a first-order
  structure for $\FOL$. Indeed, we define the interpretation of
  $R$ as $\fleq$, and the interpretation of the unary
  predicates is given via the valuations of the proposition letters.
  We write $\mo{M}^{\circ}$ for the \Lone-model $\mo{M}$ conceived of as
  a first-order structure for $\FOL$.
  The following proposition is an adaptation of~\cite[Proposition~3.26]{BezEA22}.
  Satisfaction of $\FOL$-formulae in a first-order structure is defined
  as usual.

\begin{proposition}\label{prop:st1}
  For every \Lone-model $\mo{M} = (W, 1, \fmeet, V)$, state $w \in W$ and
  formula $\phi \in \lan{L}$ we have
  $$
    \mo{M}, w \Vdash \phi \iff \mo{M}^{\circ} \models \st_x(\phi)[w].
  $$
\end{proposition}
\begin{proof}
  We use induction on the structure of $\phi$.
  If $\phi = p$ or $\phi = \top$ then the statement is obvious.
  If $\phi = \bot$ then we have
  $\mo{M}, w \Vdash \phi$ iff $w = 1$
  iff $\mo{M}^{\circ} \models \forall y(yRw)$
  iff $\mo{M}^{\circ} \models \st_x(\bot)[w]$.
  If $\phi$ is of the form $\phi_1 \wedge \phi_2$ then the inductive
  step is routine.
  If $\phi = \phi_1 \vee \phi_2$ and
  $\mo{M}, w \Vdash \phi_1 \vee \phi_2$ then there exist
  $u, v \in W$ such that $u \fmeet v \fleq w$ and
  $\mo{M}, u \Vdash \phi_1$ and $\mo{M}, v \Vdash \phi_2$.
  Taking $x' = u \fmeet v$, $y = u$ and $z = v$, this witnesses
  \begin{equation}\label{eq:st1}
    \mo{M}^{\circ} \models \exists x' \exists y \exists z(\ismeet(x', y, z)
    \wedge x'Rx \wedge st_y(\phi_1) \wedge \st_z(\phi_2))[w],
  \end{equation}
  so that $\mo{M}^{\circ} \models \st_x(\phi_1 \vee \phi_2)[w]$.
  Conversely, validity of~\eqref{eq:st1} entails the existence of
  suitable states $u$ and $v$ witnessing $\mo{M}, w \Vdash \phi_1 \vee \phi_2$.
\end{proof}

  Clearly not every structure for $\FOL$ is of the form
  $\mo{M}^{\circ}$. We can classify the ones that are.

\begin{definition}\label{def:fsl}
  Let $\FSL$ be the class of first-order structures for $\FOL$
  that satisfy the following axioms:
  \begin{enumerate}\itemindent=2em\itemsep=1pt
    \renewcommand{\labelenumi}{\textup{(\theenumi)}\enspace}
    \renewcommand{\theenumi}{M$_{\arabic{enumi}}$}
    \item \label{ax:M1}
          $\forall x (xRx)$
    \item \label{ax:M2}
          $\forall x \forall y (xRy \wedge yRx \to x = y)$
    \item \label{ax:M3}
          $\forall x \forall y \forall z(xRy \wedge yRz \to xRz)$
    \item $\forall x \forall y \exists z(zRx \wedge zRy \wedge
           \forall z'(z'Rx \wedge z'Ry \to z'Rz))$
    \item \label{it:M4b}
          $\exists x \forall y(yRx)$
    \item \label{it:M5}
          $(\exists w Pw) \wedge
          \forall x \forall y \forall z(\ismeet(x, y, z) \to ((Py \wedge Pz) \leftrightarrow Px))$
  \end{enumerate}
  Here $P$ ranges over all unary predicates of $\FOL$.
\end{definition}

  Axioms~\eqref{ax:M1} to~\eqref{ax:M3} state that the interpretation of $R$
  should be a partial order. The fourth one adds that this partial order should
  have binary greatest lower bounds and~\eqref{it:M4b} stipulates a top element.
  Finally, we have an axiom for each
  unary predicate stating that its interpretation should be a filter in the
  meet-semilattice induced by the domain and the interpretation of $R$.
  
\begin{proposition}
  A structure for $\FOL$ is isomorphic to
  $\mo{M}^{\circ}$ for some \Lone-model $\mo{M}$ iff it is in $\FSL$.
\end{proposition}
\begin{proof}
  The direction from left to right is easy.
  Conversely, let $\fmo{M} = (W, I(R), \{ I(P_p) \mid p \in \Prop \})$
  be a first-order structure that satisfies all axioms from
  Definition~\ref{def:fsl}.
  Then $I(R)$ is a partial order on $W$ with binary greatest lower bounds
  and a top element. Denote the latter by $1 \in W$ and
  the greatest lower bound of $w, v \in W$ by $w \fmeet v$.
  Define a valuation $V$ by $V(p) = I(P_p) \subseteq W$.
  Then by~\eqref{it:M5} $V$ assigns to each proposition letter a filter.
  Therefore $\fmo{M}_{\circ} := (W, 1, \fmeet, V)$ is an \Lone-model.
  It is easy to see that the identity on $W$ yields an isomorphism 
  between $\fmo{M}$ and $(\fmo{M}_{\circ})^{\circ}$.
\end{proof}
  
  Finally, we recall basic properties of $\omega$-extensions needed for the
  proof of the characterisation theorem.
  Let $\fmo{M} = (W, I(R), \{ I(P_p) \mid p \in \Prop \})$ be a $\FOL$-structure.
  For a set $A \subseteq W$, the \emph{$A$-expansion $\FOL[A]$} of $\FOL$ is
  obtained by extending $\FOL$ with new constants $\und{a}$ for each $a \in A$.
  $\fmo{M}_A$ is the expansion of $\fmo{M}$ to a structure for
  $\FOL[A]$ where $\und{a}$ is interpreted as $a \in W$.

\begin{definition}\label{def:omega-sat}
  A $\FOL$-model $\fmo{M} = (W, I(R), \{ I(P_p) \mid p \in \Prop \})$ is called
  \emph{$\omega$-saturated} if for all finite $A \subseteq W$
  and every collection $\Gamma(x_1, \ldots, x_n)$ of $\FOL[A]$-formulae with
  a finite number $n$ of free variables the following holds:
  if $\Gamma(x)$ is finitely satisfiable in $\fmo{M}_A$, then it is
  satisfiable in $\fmo{M}_A$.
  
\end{definition}

\begin{remark}
  The usual definition of $\omega$-saturation uses only one free variable.
  However, we may equivalently assume a finite number of variables,
  see e.g.~\cite[Definition~2.3.6]{ChaKei73}.
  We need the formulation with multiple variables in
  Lemma~\ref{lem:omega-mod-sat} below.
\end{remark}
  
  Using e.g.~ultraproducts, one can show that every $\FOL$-model
  has an $\omega$-saturated elementary extension~\cite{ChaKei73}.
  We denote this extension of $\fmo{M}$ by $\fmo{M}^*$,
  and the image of a state $w$ under the extension is denoted by $w^*$.
  Moreover, if $\fmo{M} \in \FSL$ then its $\omega$-saturated
  elementary extension is also in $\FSL$, since validity of the axioms
  from Definition~\ref{def:fsl} is preserved under elementary extensions.

\section{Simulations}\label{sec:4}

  In this section we define simulations between \Lone-models.
  Simulations only preserve truth of formulae in one direction, that is,
  if $S$ is a simulation and $(w, w') \in S$ then every formula satisfied
  at $w$ is also satisfied at $w'$. This prevents preservation of negations,
  and hence has been used to characterise the negation-free part of classical
  normal modal logic~\cite{KurRij97}.
  Since the language of non-distributive positive logic does not have negations,
  this is a good starting point for our attempt to characterise it.
  
  However, by its nature simulations preserve all monotone connectives,
  including classical (locally evaluated) disjunctions.
  So this approach is bound to fail,
  since the collection of first-order formulae preserved by simulations 
  is closed under classical disjunctions.
  Simulations are still worth investigating
  because they provide a stepping stone towards the Van Benthem style characterisation
  we are after.
  In fact, in Section~\ref{sec:towards-VB} we will prove that a first-order
  formula is preserved by simulations if and only if it is the classical
  disjunction of standard translations of formulae in $\lan{L}$.
  We extend this to a proper characterisation in Sections~\ref{sec:VB}
  and~\ref{sec:7b} using the notion of a meet-simulation.

\begin{definition}
  Let $\mo{M} = (W, 1, \fmeet, V)$ and $\mo{M}' = (W', 1', \fmeet', V')$ be two
  \Lone-models. An \emph{\Lone-simulation} from $\mo{M}$ to $\mo{M}'$ is a relation
  $S \subseteq W \times W'$ such that for all $(w, w') \in S$:
  \begin{enumerate}\itemindent=2em\itemsep=1pt
    \renewcommand{\labelenumi}{\textup{(\theenumi)}\enspace}
    \renewcommand{\theenumi}{S$_{\arabic{enumi}}$}
    \item \label{it:sim-1}
          If $w \in V(p)$ then $w' \in V'(p)$, for all $p \in \Prop$;
    \item \label{it:sim-2}
          If $w = 1$ then $w' = 1'$;
    \item \label{it:sim-3}
          If $v, u \in W$ are such that $v \fmeet u \fleq w$, then there
          exist $v', u' \in W'$ such that $(v, v') \in S$
          and $(u, u') \in S$ and $v' \fmeet' u' \fleq' w'$.
  \end{enumerate}
  We call $w \in W$ and $w' \in W'$ \emph{\Lone-similar} if there is an \Lone-simulation
  $S$ between $\mo{M}$ and $\mo{M}'$ such that $(w, w') \in S$.
  This is denoted by $\mo{M}, w \simul \mo{M}', w'$.
\end{definition}

  It is straightforward to see that the collection of \Lone-simulations between
  two models is closed under arbitrary unions. Therefore we have:
  
\begin{proposition}
  Let $\mo{M}$ and $\mo{M}'$ be two \Lone-models.
  \begin{enumerate}
    \item The collection of \Lone-simulations between $\mo{M}$ and $\mo{M}'$
          forms a complete lattice.
    \item The relation of \Lone-similarity between $\mo{M}$ and $\mo{M}'$ is again
          an \Lone-simulation.
  \end{enumerate}
\end{proposition}

  There are several generic examples of \Lone-simulations.

\begin{example}\label{exm:easy}
  Let $\mo{M} = (W, 1, \fmeet, V)$ be an \Lone-model.
  Then the relations 
  $S_1 = \{ (w, w) \mid w \in W \}$ and
  $S_2 = \{ (w, v) \in W \mid w \fleq v \}$
  are \Lone-simulations from $\mo{M}$ to itself.
\end{example}

  \Lone-morphisms between models also give rise to \Lone-simulations.

\begin{example}
  Let $\mo{M} = (W, 1, \fmeet, V)$ and $\mo{M}' = (W', 1', \fmeet', V')$
  be \Lone-models and $f : \mo{M} \to \mo{M}'$ an \Lone-morphism.
  Then the following variations on the graph of $f$ are \Lone-simulations:
  \begin{align*}
    \Gr f &:= \{ (w, f(w)) \in W \times W' \mid w \in W \}, \\
    \Gr^{\uparrow} f &:= \{ (w, w') \in W \times W' \mid f(w) \fleq' w' \}.
  \end{align*}
  Taking $f = \id : \mo{M} \to \mo{M}$ we obtain $S_1$ and $S_2$
  from Example~\ref{exm:easy} as $\Gr f$ and $\Gr^{\uparrow}f$.
\end{example}

  As expected, \Lone-similarity implies logical inclusion.
  
\begin{proposition}\label{prop:adeq}
  If $\mo{M}, w \simul \mo{M}', w'$ then $\mo{M}, w \login \mo{M}', w'$.
\end{proposition}
\begin{proof}
  Let $\phi \in \lan{L}$ be such that $\mo{M}, w \Vdash \phi$.
  We show by induction on the structure of $\phi$ that $\mo{M}', w' \Vdash \phi$.
  The case for $\phi = \top$ is trivial, and
  if $\phi = p \in \Prop$ or $\phi = \bot$ then the result follows
  from~\eqref{it:sim-1} and~\eqref{it:sim-2}, respectively. 
  The inductive step for $\phi = \phi_1 \wedge \phi_2$ is straightforward.

  Suppose $\mo{M}, w \Vdash \phi_1 \vee \phi_2$.
  Then there are $u, v \in W$ such that
  $\mo{M}, u \Vdash \phi_1$ and $\mo{M}, v \Vdash \phi_2$ and
  $u \fmeet v \fleq w$.
  Using~\eqref{it:sim-3} we can find
  $u', v' \in W'$ such that $(u, u') \in S$ and $(v, v') \in S$ and
  $u' \fmeet' v' \fleq' w'$.
  By the induction hypothesis $\mo{M}', u' \Vdash \phi_1$ and
  $\mo{M}', v' \Vdash \phi_2$, which implies
  $\mo{M}', w' \Vdash \phi_1 \vee \phi_2$.
\end{proof}

\section{A Hennessy-Milner style theorem}\label{sec:5}

  While Proposition~\ref{prop:adeq} guarantees that \Lone-similarity implies
  logical inclusion, the opposite (i.e.~logical inclusion implies \Lone-similarity)
  need not be true.
  Classes on which \Lone-similarity and logical inclusion coincide are often
  called Hennessy-Milner classes, after the authors who proved an
  analogous result for classical normal modal logic~\cite{HenMil85}.
  This section is devoted to finding Hennessy-Milner classes.

\begin{definition}
  A class $\ms{C}$ of \Lone-models is called a \emph{Hennessy-Milner class}
  if for all $\mo{M}, \mo{M}' \in \ms{C}$ and states $w, w'$ in them
  we have
  $$
    \mo{M}, w \login \mo{M}', w'
      \iff \mo{M}, w \simul \mo{M}', w'.
  $$
\end{definition}

  We start by proving that the class of finite models is a
  Hennessy-Milner class. This will serve as inspiration for an analogue of
  modal saturation which we formulate in terms of compactness
  in a topology, and a more general Hennessy-Milner style result.

\begin{proposition}\label{prop:HM-fin}
  The class of finite \Lone-models is a Hennessy-Milner class.
\end{proposition}
\begin{proof}
  Let $\mo{M} = (W, 1, \fmeet, V)$ and $\mo{M}' = (W', 1', \fmeet', V')$ be two
  finite \Lone-models. We claim that
  $$
    S := \{ (w, w') \in W \times W' \mid \mo{M}, w \login \mo{M}', w' \}
  $$
  is an \Lone-simulation. Together with Proposition~\ref{prop:adeq}
  this proves the proposition.
  Item~\eqref{it:sim-1} holds by definition of $S$,
  and~\eqref{it:sim-2} follows from the fact that $1$ and $1'$ are the
  only elements satisfying $\bot$.
  
  Let $(w, w') \in S$ and let $v, u \in W$ be such that $v \fmeet u \fleq w$.
  Suppose towards a contradiction that~\eqref{it:sim-3} is not satisfied.
  Then for each
  $(v', u') \in M(w') := \{ (v', u') \in W' \times W' \mid v' \fmeet' u' \fleq' w' \}$
  either
  \begin{itemize}
    \item $(v, v') \notin S$, so there exists a formula $\phi$
          that is satisfied at $v$ but not at $v'$; or
    \item $(u, u') \notin S$, so there exists a formula $\psi$
          that is satisfied at $u$ but not at $u'$.
  \end{itemize}
  For each $(v', u') \in M(w')$ select a $\phi$ or $\psi$ as specified.
  Let $\Phi$ be the set of all such $\phi$ and $\Psi$ of all such $\psi$.
  Let $\bar{\phi} := \bigwedge \Phi$
  and $\bar{\psi} := \bigwedge \Psi$
  (taking to empty conjunction to be $\top$).
  Then $v$ and $u$ satisfy $\bar{\phi}$ and $\bar{\psi}$, respectively.
  On the other hand, for each $(v', u') \in M(w')$ either $v'$
  does not satisfy $\bar{\phi}$ or $u'$ does not satisfy $\bar{\psi}$.
  Therefore
  $$
    \mo{M}, w \Vdash \bar{\phi} \vee \bar{\psi}
    \quad\text{ but }\quad
    \mo{M}', w' \not\Vdash \bar{\phi} \vee \bar{\psi}.
  $$
  This contradicts the assumption that $(w, w') \in S$.
  So~\eqref{it:sim-3} must hold and $S$ is an \Lone-simulation.
\end{proof}

  If we try to apply this proof to infinite models we may encounter
  infinite sets $\Phi$ and $\Psi$. So our analogue of modal saturation
  should remedy this.
  Indeed, we need a compactness property that allows us to reduce
  infinite $\Phi$ and $\Psi$ to finite sets.

\begin{definition}
  Let $\mo{M} = (W, 1, \fmeet, V)$ be an \Lone-model.
  Then we denote by $\tau_V$ the topology on $W$ generated by
  the clopen subbase
  $$
    \{ \llb \phi \rrb^{\mo{M}} \mid \phi \in \lan{L} \}
    \cup
    \{ W \setminus \llb \phi \rrb^{\mo{M}} \mid \phi \in \lan{L} \}.
  $$
  The \Lone-model $\mo{M}$ is called \emph{meet-compact} if for
  all $w \in W$ the set
  $
    M(w) := \{ (v, u) \in W \times W \mid v \fmeet u \fleq w \}
  $
  is a compact subset of $(W, \tau_V) \times (W, \tau_V)$.
\end{definition}

\begin{example}
  Every finite monotone \Lone-model is meet-compact.
  Evidently, the finiteness entails the compactness requirement.
\end{example}
 
  As another example, we show that $\omega$-saturated \Lone-models
  are meet-compact. An \Lone-model $\mo{M}$ is called $\omega$-saturated
  if $\mo{M}^{\circ}$ is $\omega$-saturated (see~Definition~\ref{def:omega-sat}).

\begin{lemma}\label{lem:omega-mod-sat}
  If $\mo{M} = (W, 1, \fmeet, V)$ is $\omega$-saturated,
  then it is meet-compact.
\end{lemma}
\begin{proof}
  Let $w \in W$. We need to show that
  $M(w) = \{ (u, v) \in W \times W \mid u \fmeet v \fleq w \}$ is compact
  in $(W, \tau_V) \times (W, \tau_V)$.
  By the Alexander subbase theorem, it suffices to show that every
  open cover of subbasic opens has a finite subcover.
  A subbase for the topology on $(X, \tau_V) \times (X, \tau_V)$ can be
  given by the collection of open squares $a \times b$, where $a$ and $b$
  range over the subbase for $\tau_V$, that is, over the truth sets of
  formulae and their complements.
  So suppose
  \begin{equation}\label{eq:open-cover}
  \begin{split}
    M(w) \subseteq &\bigcup \{ \llb \phi_i \rrb^{\mo{M}} \times 
           \llb \psi_i \rrb^{\mo{M}} \mid i \in I \}
      \cup \bigcup \{ \llb \phi_j \rrb^{\mo{M}} \times
           (X \setminus \llb \psi_j \rrb^{\mo{M}}) \mid j \in J \} \\
      &\cup \bigcup \{ (X \setminus \llb \phi_k \rrb^{\mo{M}}) \times
           \llb \psi_k \rrb^{\mo{M}} \mid k \in K \}
      \cup \bigcup \{ (X \setminus \llb \phi_{\ell} \rrb^{\mo{M}}) \times
           (X \setminus \llb \psi_{\ell} \rrb^{\mo{M}}) \mid \ell \in L \}
  \end{split}
  \end{equation}
  where $I, J, K, L$ are index sets.
  Let $A = \{ \und{w} \}$ and
  \begin{align*}
    \Gamma(x, y, z) = \{ x R \und{w} \}
      & \cup \{ \ismeet(x, y, z) \} \\
      &\cup \bigcup \big\{ 
        \neg\st_y(\phi_i) \wedge \neg\st_z(\psi_i) \mid i \in I \big\} \\
      &\cup \bigcup \big\{
        \neg\st_y(\phi_j) \wedge \st_z(\psi_j) \mid j \in J \big\} \\
      &\cup \bigcup \big\{
        \st_y(\phi_k) \wedge \neg\st_z(\psi_k) \mid k \in K \big\} \\
      &\cup \bigcup \big\{
        \st_y(\phi_{\ell}) \wedge \st_z(\psi_{\ell}) \mid \ell \in L \big\}.
  \end{align*}
  Then $\Gamma(x)$ is not satisfiable in $\mo{M}^{\circ}_A$,
  because if it were then there is a state $x \in W$ such that $x R w$ which is the meet
  of a pair $(y, z)$ that is not in the open cover given in
  \eqref{eq:open-cover}.
  Since $\mo{M}^{\circ}$ is $\omega$-saturated, this implies that $\Gamma(x)$ is not
  finitely satisfiable in $\mo{M}^{\circ}_A$,
  which by a reverse argument yields a finite subcover of~\eqref{eq:open-cover}.
\end{proof}

\begin{theorem}\label{thm:HM}
  The meet-compact \Lone-models form a Hennessy-Milner class.
\end{theorem}
\begin{proof}
  Let $\mo{M} = (W, 1, \fmeet, V)$ and $\mo{M}' = (W', 1', \fmeet', V')$ be two
  meet-compact \Lone-models.
  We claim that
  $$
    S = \{ (w, w') \in W \times W' \mid \mo{M}, w \login \mo{M}', w' \}
  $$
  is an \Lone-simulation. Together with Proposition~\ref{prop:adeq} this proves
  the theorem.
  
  Item~\eqref{it:sim-1} is satisfied by definition,
  and~\eqref{it:sim-2} follows from the fact that $1$ and $1'$ are the
  only elements satisfying $\bot$.
  To prove~\eqref{it:sim-3}, assume $(w, w') \in S$
  and $u, v \in W$ are such that $u \fmeet v \fleq w$.
  Suppose there are no $u', v' \in W'$ such that
  $u' \fmeet' v' \fleq' w'$ and $(u, u') \in S$ and $(v, v') \in S$.
  Proceeding as in the proof of Proposition~\ref{prop:HM-fin},
  we obtain (potentially infinite) sets $\Phi$ and $\Psi$.
  Now observe that we have an open cover
  $$
    M(w')
      \subseteq \bigcup_{\phi \in \Phi}
          \big((W' \setminus \llb \phi \rrb^{\mo{M}'}) \times W' \big)
        \cup \bigcup_{\psi \in \Psi}
          \big(W' \times (W' \setminus \llb \psi \rrb^{\mo{M}'})\big).
  $$
  By assumption $M(w')$ is compact,
  so we can find a finite subcover indexed by finite sets
  $\Phi' \subseteq \Phi$ and $\Psi' \subseteq \Psi$.
  Define $\bar{\phi} := \bigwedge \Phi'$ and
  $\bar{\psi} := \bigwedge \Psi'$.
  The remainder of the proof is the analogous to
  Proposition~\ref{prop:HM-fin}.
\end{proof}

\section{Towards characterisation}\label{sec:towards-VB}

  In this section we characterise the \Lone-simulation-invariant fragment of
  $\FOL$ over the class $\FSL$ of first-order structures corresponding
  to \Lone-models.

\begin{definition}
  A first-order formula $\alpha(x)$ with one free variable $x$
  is \emph{preserved by \Lone-simulations}
  if 
  $$
    \mo{M}^{\circ} \models \alpha(x)[w]
    \quad\text{implies}\quad
    \mo{N}^{\circ} \models \alpha(x)[v]
  $$
  whenever $\mo{M}, w \simul \mo{N}, v$, for all \Lone-models $\mo{M}, \mo{N}$.
\end{definition}

  We cannot yet characterise the \Lone-simulation-preserving formulae as the
  language $\lan{L}$, because we can find first-order formulae $\alpha(x)$
  that are preserved by \Lone-simulations but not equivalent to the
  standard translation of a formula in $\lan{L}$.
  The next example gives such an $\alpha(x)$.

\begin{example}\label{exm:6}
  Let $p$ and $q$ be proposition letters, and $P$ and $Q$ their
  corresponding unary predicates. Consider the first-order
  formula $\alpha(x) := Px \cvee Qx$
  (recall that $\cvee$ denotes classical disjunction).
  If $\mo{M}, w \simul \mo{N}, v$ and
  $\mo{M}^{\circ} \models \alpha(x)[w]$, then
  $\mo{M}, w \Vdash p$ or $\mo{M}, w \Vdash q$, which by
  \eqref{it:sim-1} implies that $\mo{N}, v \Vdash p$ or $\mo{N}, v \Vdash q$.
  Hence $\mo{N} \models \alpha(x)[v]$. So $\alpha(x)$ is preserved by
  \Lone-simulations.
  
  We give a model $\mo{M}$ where the collection of states $w$
  such that $\mo{M}^{\circ} \models \alpha(x)[w]$ is not a filter.
  Since the interpretation of any formula $\phi \in \lan{L}$ yields
  a filter, it then follows from Proposition~\ref{prop:st1} that
  $\alpha(x)$ is not equivalent to the standard translation of a formula
  in $\lan{L}$.
  Consider the model $\mo{M} = (W, 1, \fmeet, V)$,
  where $W = \{ w, u, v, 1 \}$ is ordered by $w \fleq u \fleq 1$
  and $w \fleq v \fleq 1$.
  Let $V(p) = \{ u, 1 \}$ and $V(q) = \{ v, 1 \}$.
  Then $\mo{M}$ is an \Lone-model 
  and the set $\{ w \in W \mid \mo{M}^{\circ} \models \alpha(x)[w] \}
  = \{ u, v, 1 \}$
  is not a filter in $(W, 1, \fmeet)$ (see Figure~\ref{fig:exm-trisim-cjoin}).
  \begin{figure}[H]
    \centering
    \begin{tikzpicture}[scale=.9]
      \draw[rotate=-45] (-1.7,.7) ellipse(.7 and 1.3);
      \draw (-1,1.9) node{$p$};
      \draw[rotate=45] (1.7,.7) ellipse(.7 and 1.3);
      \draw (1,1.9) node{$q$};
      \draw[fill=black] (0,0) circle(.5mm);
      \draw (0,-.25) node{$w$};
      \draw[fill=black] (-1,1) circle(.5mm);
      \draw (-1.25,1) node{$u$};
      \draw[fill=black] (1,1) circle(.5mm);
      \draw (1.25,1) node{$v$};
      \draw[fill=black] (0,2) circle(.5mm);
      \draw (0,2.25) node{$1$};
      \draw[-latex,thick] (.11,.11) -- (.89,.89);
      \draw[-latex,thick] (-.11,.11) -- (-.89,.89);
      \draw[-latex,thick] (-.89,1.11) -- (-.11,1.89);
      \draw[-latex,thick] (.89,1.11) -- (.11,1.89);
    \end{tikzpicture}
    \caption{The model $\mo{M} = (W, 1, \fmeet, V)$ from Example~\ref{exm:6}.}
    \label{fig:exm-trisim-cjoin}
  \end{figure}
\end{example}

  Along the lines of~\cite[Theorem~2.68]{BRV01}, we characterise the first-order formulae
  preserved by simulations as those formulae in $\FOL$
  equivalent to the disjunction of standard translations of $\lan{L}$-formulae.

\begin{theorem}\label{thm:towards-VB}
  A first-order formula $\alpha(x) \in \lan{FOL}$ with one free variable $x$
  is equivalent over $\cat{FSL}$
  to a formula of the form $\st_x(\phi_1) \cvee \cdots \cvee \st_x(\phi_n)$,
  where $\phi_i \in \lan{L}$, if and only if it is
  preserved by \Lone-simulations.
\end{theorem}
\begin{proof}
  The left-to-right implication follows from Proposition~\ref{prop:adeq}.
  For the converse, let $\alpha(x)$ be a formula preserved by \Lone-simulations.
  For $\phi_1, \ldots, \phi_n \in \lan{L}$, abbreviate
  $\st_x[\phi_1, \ldots, \phi_n] := \st_x(\phi_1) \cvee \cdots \cvee \st_x(\phi_n)$.
  (We take the empty disjunction to be the formula $\neg(x = x)$ in $\FOL$.)
  Define
  $$
    \MOC(\alpha)
      = \{ \st_x[\phi_1, \ldots, \phi_n]
           \mid \alpha(x) \models_{\cat{FSL}} \st_x[\phi_1, \ldots, \phi_n]
           \text{ and } n \in \omega \text{ and } \phi_i \in \lan{L} \}.
  $$
  It suffices to prove that $\MOC(\alpha) \models_{\cat{FSL}} \alpha(x)$, because
  a compactness argument then yields the existence of a finite subset
  $\MOC'(\alpha) := \{ \st_x[\phi_{1,1}, \ldots, \phi_{1,n_1}], \ldots, \st_x[\phi_{k,1}, \ldots, \phi_{k,n_k}] \} \subseteq \MOC(\alpha)$
  that entails $\alpha(x)$.
  Since also $\alpha \models_{\cat{FSL}} \st_x[\phi_{i,1}, \ldots, \phi_{i,n_i}]$
  by definition, we then find that $\alpha(x)$ is equivalent over $\cat{FSL}$ to
  \begin{align*}
    \st_x&[\phi_{1,1}, \ldots, \phi_{1,n_1}]
      \wedge \cdots \wedge \st_x[\phi_{k,1}, \ldots, \phi_{k,n_k}] \\
      &= \bigvee^{c} \big\{ \st_x(\phi_{1, v(1)}) \wedge \cdots \wedge \st_x(\phi_{k, v(k)})
         \mid v(i) \in \{ 1, \ldots, n_i \}
              \text{ for all } i \in \{ 1, \ldots, k \} \big\} \\
      &= \bigvee^c \big\{ \st_x(\phi_{1, v(1)} \wedge \cdots \wedge \phi_{k, v(k)})
         \mid v(i) \in \{ 1, \ldots, n_i \}
              \text{ for all } i \in \{ 1, \ldots, k \} \big\}
  \end{align*}
  For the first equality above we use distributivity in $\FOL$.
  The resulting formula is the finite disjunction of standard translations
  of formulae in $\lan{L}$, because
  $\phi_{1, v(1)} \wedge \cdots \wedge \phi_{k, v(k)} \in \lan{L}$
  whenever $\phi_{1, v(1)}, \ldots, \phi_{k, v(k)} \in \lan{L}$.
  
  So let $\fmo{M} \in \cat{FSL}$ and assume $\fmo{M} \models \MOC(\alpha)[w]$.
  Let $\fmo{M}_{\circ} = (W, 1, \fmeet, V)$. 
  For $w \in W$, let
  \begin{equation}\label{eq:satisf}
    \neg\th(w)
      := \{ \phi \in \lan{L} \mid
            \fmo{M}_{\circ}, w \not\Vdash \phi \}.
  \end{equation}
  We claim that $\{ \alpha(x) \} \cup \{ \neg\st_x(\phi) \mid \phi \in \neg\th(w) \}$
  is satisfiable.
  Suppose not, then there exists a finite number of formulae
  $\phi_1, \ldots, \phi_m \in \neg\th(w)$ such that
  $\alpha(x) \models \neg(\neg\st_x(\phi_1) \wedge \cdots \wedge \neg\st_x(\phi_m))$.
  This implies $\alpha(x) \models \st_x(\phi_1) \cvee \cdots \cvee \st_x(\phi_m)$.
  But then $\st_x(\phi_1) \cvee \cdots \cvee \st_x(\phi_m) \in \MOC(\alpha)$,
  so by assumption $\fmo{M} \models \st_x(\phi_1) \cvee \cdots \cvee \st_x(\phi_m)[w]$.
  This implies that $\fmo{M} \models \st_x(\phi_i)[w]$ for one of the $\phi_i$,
  and hence $\fmo{M}_{\circ}, w \Vdash \phi_i$, which contradicts
  the assumption that $\phi_i \in \neg\th(w)$.
  So the set in \eqref{eq:satisf} is satisfiable.
  
  Thus, we can find a structure $\fmo{N} \in \cat{FSL}$ and a state
  $v$ in the domain of $\fmo{N}$ such that $\fmo{N} \models \alpha(x)[v]$ and
  $\fmo{N}_{\circ}, v_{\circ} \not\Vdash \phi$ for all $\phi \in \neg\th(w)$.
  This implies
  \begin{equation}\label{eq:VB-login}
    \fmo{N}_{\circ}, v \login \fmo{M}_{\circ}, w.
  \end{equation}
  Recall that we write $\fmo{M}^*$ and $\fmo{N}^*$ for $\omega$-saturated
  elementary extensions of $\fmo{M}$ and $\fmo{N}$,
  and $w^*$ and $v^*$ for the images of $w$ and $v$.
  As a consequence of Lemma~\ref{lem:omega-mod-sat} the models $(\fmo{M}^*)_{\circ}$ and
  $(\fmo{N}^*)_{\circ}$ are meet-compact. Since elementary embeddings
  preserve truth of formulae, \eqref{eq:VB-login} and Theorem~\ref{thm:HM}
  imply
  $$
    (\fmo{N}^*)_{\circ}, v^* \simul (\fmo{M}^*)_{\circ}, w^*.
  $$
  Now $\fmo{N} \models \alpha(x)[v]$ implies $\fmo{N}^* \models \alpha(x)[v^*]$.
  Since $\alpha(x)$ is assumed to be preserved by \Lone-simulations,
  we find $\fmo{M}^* \models \alpha(x)[w^*]$.
  Invariance of truth of formulae under elementary embeddings
  gives $\fmo{M} \models \alpha(x)[w]$.
\end{proof}

\begin{example}\label{exm:sim-VB-join}
  We investigate what goes wrong in the proof of Theorem~\ref{thm:towards-VB}
  if we try to eliminate the additional classical disjunctions by taking
  $$
    \MOC(\alpha)
      = \{ \st_x(\phi) \mid \alpha(x) \models_{\cat{FSL}} \st_x(\phi), \phi \in \lan{L} \}.
  $$
  Let $\alpha(x) = Px \cvee Qx$.
  Since $p \vee q$ is the smallest $\lan{L}$-formula whose truth
  set contains both $\llb p \rrb$ and $\llb q \rrb$ (in any model),
  any formula $\phi \in \lan{L}$ such that $\alpha(x) \models_{\FSL} \st_x(\phi)$
  must be implied by $p \vee q$ (in the sense that truth of $p \vee q$ implies
  truth of $\phi$).
  Taking $\mo{M}$ as in Example~\ref{exm:6}, we find
  $\mo{M}, w \Vdash p \vee q$ hence $\mo{M}^{\circ} \models \st_x(p \vee q)[w]$.
  Since we argued that $\st_x(p \vee q)$ is the smallest formula
  in $\MOC(\alpha)$, we have $\mo{M}^{\circ} \models \MOC(\alpha)[x]$.
  But at the same time both $p$ and $q$ are in $\neg\th(w)$.
  Then supposing that $\{ \alpha(x) \} \cup \{ \neg\st_x(\phi) \mid \phi \in \neg\th(x) \}$
  is not satisfiable, we obtain $\alpha \models_{\cat{FSL}} \st_x(p) \cvee \st_x(q)$.
  This does not yield a contradiction, because $\st_x(p) \cvee \st_x(q)$ is not
  in $\MOC(\alpha)$.
  Hence the proof fails.
\end{example}

\section{Characterisation using meet-simulations}\label{sec:VB}

  Making use of Theorem~\ref{thm:towards-VB}, we now prove a
  Van Benthem style characterisation theorem for $\lan{L}$.
  As we have seen, the main challenge is to prevent classical disjunctions
  from being preserved. We do so by adapting the notion of an \Lone-simulation to
  that of a meet-simulation.
  Instead of relating states with states, meet-simulations relate
  pairs of states from one model to single states of another.

\begin{definition}\label{def:meet-sim}
  Let $\mo{M} = (W, 1, \fmeet, V)$ and $\mo{M}' = (W', 1', \fmeet', V')$
  be two \Lone-models.
  A \emph{meet-simulation} from $\mo{M}$ to $\mo{M}'$
  is a relation $T \subseteq (W \times W) \times W'$ such that for
  all $(w_1, w_2, w') \in T$:
  \begin{enumerate}\itemindent=2em\itemsep=1pt
    \renewcommand{\labelenumi}{\textup{(\theenumi)}\enspace}
    \renewcommand{\theenumi}{M$_{\arabic{enumi}}$}
    \item \label{it:msim-1}
          If $w_1 \in V(p)$ and $w_2 \in V(p)$ then $w' \in V'(p)$, for all $p \in \Prop$;
    \item \label{it:msim-2}
          If $w_1 = w_2 = 1$ then $w' = 1'$;
    \item \label{it:msim-3}
          If $u_1, v_1, u_2, v_2 \in W$ are such that $u_1 \fmeet v_1 \fleq w_1$
          and $u_2 \fmeet v_2 \fleq w_2$, then there exist $v', u' \in W'$ such that
                  $(u_1, u_2, u') \in T$ and $(v_1, v_2, v') \in T$
                  and $v' \fmeet' u' \fleq' w'$.
  \end{enumerate}
\end{definition}

\begin{remark}
  Instead of relating pairs of states from one model to states of another,
  we can also define meet-simulations as relations between three (possibly
  distinct) models. The definition above is then obtained as the special case
  where the first two models are the same.
  The results in this section work for either definition of meet-simulation.
\end{remark}

  We give some examples of meet-simulations.

\begin{example}\label{exm:meet-sim-on-M}
  For any \Lone-model $\mo{M} = (W, 1, \fmeet, V)$, the relation
  $$
    T = \{ (w_1, w_2, w_3) \in W \times W \times W
           \mid w_1 \fmeet w_2 \fleq w_3 \}
  $$
  is a meet-simulation on $\mo{M}$.
  Let us verify this.
  Condition~\eqref{it:msim-1} follows from the fact that proposition
  letters are interpreted as filters, so if both $w_1$ and $w_2$ satisfy
  $p \in \Prop$, then so does $w_1 \fmeet w_2$ and everything above it.
  Condition~\eqref{it:msim-2} follows immediately from the definition.
  For~\eqref{it:msim-3}, suppose $(w_1, w_2, w_3) \in T$,
  $u_1 \fmeet v_1 \fleq w_1$ and $u_2 \fmeet v_2 \fleq w_2$.
  Then $u_3 := u_1 \fmeet u_2$ and $v_3 := v_1 \fmeet v_2$
  witness truth of~\eqref{it:msim-3}, as
  $u_3 \fmeet v_3
    = u_1 \fmeet u_2 \fmeet v_1 \fmeet v_2
    \fleq w_1 \fmeet w_2 \fleq w_3$.
\end{example}

\begin{example}
  Let $\mo{M} = (W, 1, \fmeet, V)$ and $\mo{M}' = (W', 1', \fmeet', V')$
  be two \Lone-models and $f : \mo{M} \to \mo{M}'$ an \Lone-morphism.
  Then the following relations are meet-simulations:
  \begin{align*}
    T_1 &= \{ (w_1, w_2, w') \in W \times W \times W' \mid f(w_1 \fmeet w_2) = w' \} \\
    T_2 &= \{ (w_1, w_2, w') \in W \times W \times W' \mid f(w_1 \fmeet w_2) \fleq' w' \} \\
    T_3 &= \{ (w_1', w_2', w) \in W' \times W' \times W
              \mid w_1' \fmeet' w_2' \fleq' f(w) \}
  \end{align*}
  We show that $T_2$ is a meet-simulation.
  The verification for the other two is similar.
  Let $(w_1, w_2, w') \in T_2$.
  \begin{enumerate}\itemsep=1pt
    \item[\eqref{it:msim-1}]
          Suppose $w_1, w_2 \in V(p)$ for some proposition letter $p$.
          Then $w_1 \fmeet w_2 \in V(p)$ because $V(p)$ is a filter of
          $(W, 1, \fmeet)$. By definition of an \Lone-morphism this implies that
          $f(w_1 \fmeet w_2) \in V'(p)$. Using the fact that $V'(p)$ is a
          filter of $(W', 1', \fmeet')$ we find $w' \in V'(p)$.
    \item[\eqref{it:msim-2}]
          If $w_1 = w_2 = 1$ then $1' = f(1 \fmeet 1) \fleq w'$,
          which implies $w' = 1'$.
    \item[\eqref{it:msim-3}]
          Suppose that $u_1, v_1, u_2, v_2 \in W$
          are such that $u_1 \fmeet v_1 \fleq w_1$ and $u_2 \fmeet v_2 \fleq w_2$.
          Set $u' := f(u_1 \fmeet u_2)$ and $v' := f(v_1 \fmeet v_2)$.
          Then we have $(u_1, u_2, u') \in T_2$ and
          $(v_1, v_2, v') \in T_2$ by definition, and
          $$
            u' \fmeet' v' = f(u_1 \fmeet u_2) \fmeet' f(v_1 \fmeet v_2)
              = f(u_1 \fmeet v_1 \fmeet u_2 \fmeet v_2)
              \fleq' f(w_1 \fmeet w_2)
              \fleq' w'.
          $$
  \end{enumerate}
\end{example}

  We define preservation of a first-order formula with one free variable
  by meet-simulations as follows.

\begin{definition}
  A first-order formula $\alpha(x)$ with one free variable $x$
  is said to be preserved by meet-simulations if for every meet-simulation
  $T$ between $\mo{M}$ and $\mo{M}'$ with $(w_1, w_2, w') \in T$
  we have:
  $$
    \text{if}\quad \mo{M}^{\circ} \models \alpha(x)[w_1]
    \quad\text{and}\quad \mo{M}^{\circ} \models \alpha(x)[w_2]
    \quad\text{then}\quad (\mo{M}')^{\circ} \models \alpha(x)[w'].
  $$
\end{definition}

\begin{example}\label{exm:trisim-cjoin}
  Let $\mo{M} = (W, 1, \fmeet, V)$ be the model from Example~\ref{exm:6}
  and $\alpha(x) = Px \cvee Qx$.
  By Example~\ref{exm:meet-sim-on-M} the tuple
  $(u, v, w)$ is related by a meet-simulation.
  We have seen that $\mo{M}^{\circ} \models \alpha(x)[u]$
  and $\mo{M}^{\circ} \models \alpha(x)[v]$.
  But clearly we do not have $\mo{M}^{\circ} \models \alpha(x)[w]$.
  So $\alpha(x)$ is not preserved by meet-simulations.
\end{example}

\begin{proposition}\label{prop:adeq-meet-sim}
  Let $T$ be a meet-simulation between $\mo{M}$ and $\mo{M}'$,
  $(w_1, w_2, w') \in T$ and $\phi \in \lan{L}$.
  $$
    \text{If}\quad
    \mo{M}, w_1 \Vdash \phi
    \quad\text{and}\quad
    \mo{M}, w_2 \Vdash \phi
    \quad\text{then}\quad
    \mo{M}', w' \Vdash \phi.
  $$
\end{proposition}
\begin{proof}
  The proof by induction on the structure of $\phi$
  is analogous to the proof of Proposition~\ref{prop:adeq}.
\end{proof}

  In the remainder of this section we work towards the desired
  characterisation theorem.
  The following proposition allows us to exploit the result from
  Section~\ref{sec:towards-VB}.

\begin{proposition}\label{prop:sim-trisim}
  Let $\mo{M} = (W, 1, \fmeet, V)$ and $\mo{M}' = (W', 1', \fmeet', V')$
  be \Lone-models, $S \subseteq W \times W'$ and
  $$
    T_S = \{ (w_1, w_2, w') \in W \times W \times W'
           \mid (w_1, w') \in S \text{ or } (w_2, w') \in S \}.
  $$
  Then $S$ is an \Lone-simulation between $\mo{M}$ and $\mo{M}'$ if and only if
  $T_S$ is a meet-simulation between them.
\end{proposition}
\begin{proof}
  Suppose $S$ is an \Lone-simulation. We verify that $T_S$ is a meet-simulation.
  Let $(w_1, w_2, w') \in T_S$ and assume without loss of generality that
  $(w_1, w') \in S$.
  If $w_1 \in V(p)$ and $w_2 \in V(p)$
  then~\eqref{it:sim-1} implies that $w' \in V'(p)$,
  so that~\eqref{it:msim-1} is satisfied.
  If $(1, 1, w') \in T_S$ then $(1, w') \in S$ so by~\eqref{it:sim-2}
  $w' = 1'$, proving~\eqref{it:msim-2}.
  Lastly, if $v_1 \fmeet u_1 \fleq w_1$ and $v_2 \fmeet u_2 \fleq w_2$
  then using~\eqref{it:sim-3} we can find $v', u' \in W'$ such
  that $v' \fmeet' u' \fleq' w'$ and $(v_1, v') \in S$ and $(u_1, u') \in S$.
  This implies $(v_1, v_2, v') \in T_S$ and
  $(u_1, u_2, u') \in T_S$, so that~\eqref{it:msim-3} holds.

  For the reverse direction, assume that $T_S$ is a meet-simulation.
  Using the fact that for each $(w, w') \in S$ we have $(w, w, w') \in T_S$
  we can show that~%
  \eqref{it:msim-1} implies~\eqref{it:sim-1},
  \eqref{it:msim-2} implies~\eqref{it:sim-2},
  and \eqref{it:msim-3} implies~\eqref{it:sim-3}.
\end{proof}

\begin{corollary}
  If $\alpha(x)$ is preserved by meet-simulations, then it is also
  preserved by simulations.
\end{corollary}
\begin{proof}
  Let $S$ be a simulation between $\mo{M}$ and $\mo{M}'$
  such that $(w, w') \in S$ and $\mo{M}^{\circ} \models \alpha(x)[w]$.
  Constructing $T_S$ from $S$ as in Proposition~\ref{prop:sim-trisim} we
  find a meet-simulation $T_S$ such that $(w, w, w') \in T_S$.
  The assumption that $\alpha(x)$ is preserved by meet-simulations
  now entails $\mo{M}' \models \alpha(w)[w']$.
\end{proof}

\begin{lemma}\label{lem:alphax-filter}
  Let $\alpha(x)$ be a first-order formula with one free variable that is
  preserved by meet-simulations.
  Let $\fmo{M} = (W, I(R), \{ I(P_p) \mid p \in \Prop \})$ be a
  first order structure in $\FSL$ and $\fmo{M}_{\circ} = (W, 1, \fmeet, V)$.
  Then the set
  $\llb \alpha(x) \rrb^{\fmo{M}} := \{ w \in W \mid \fmo{M} \models \alpha(x)[w] \}$
  is either empty or a filter of $(W, 1, \fmeet)$. 
\end{lemma}
\begin{proof}
  If there are no $w \in W$ such that $\fmo{M} \models \alpha(x)[w]$ then
  $\llb \alpha(x) \rrb^{\fmo{M}}$ is empty.
  Assume $\llb \alpha(x) \rrb^{\fmo{M}} \neq \emptyset$.
  Define $T = \{ (w, v, u) \in W \mid (w \fmeet v) \fleq u \}$.
  We have seen in Example~\ref{exm:meet-sim-on-M} that this is a meet-simulation.
  Suppose $\fmo{M} \models \alpha(x)[w]$ and $(w, v) \in R$.
  Then $(w, w, v) \in T$ and since $\alpha(x)$ is assumed to be preserved
  by meet-simulations we find $\fmo{M} \models \alpha(x)[v]$.
  This implies that $\llb \alpha(x) \rrb^{\fmo{M}}$
  is upward closed under the partial order $\fleq$ ($= I(R)$).
  Next suppose $\fmo{M} \models \alpha(x)[w]$ and $\fmo{M} \models \alpha(x)[v]$.
  Then $(w, v, w \fmeet v) \in T$ so $\fmo{M} \models \alpha(x)[w \fmeet v]$.
  This implies that $\llb \alpha(x) \rrb^{\fmo{M}}$ is closed under binary
  meets. So $\llb \alpha(x) \rrb^{\fmo{M}}$ is a filter
  of $(W, 1, \fmeet)$.
\end{proof}

\begin{theorem}\label{thm:VB}
  Let $\alpha(x) \in \FOL$ be a formula with one free variable $x$.
  Then $\alpha(x)$ is equivalent over $\cat{FSL}$ to $\neg(x \neq x)$ or
  to the standard translation of an $\lan{L}$-formula iff it is
  preserved by meet-simulations.
\end{theorem}
\begin{proof}
  Since $\alpha(x)$ is preserved by meet-simulations, it is also 
  preserved by simulations.
  Therefore Theorem~\ref{thm:towards-VB} implies that we can find
  $\phi_1, \ldots, \phi_n$ such that
  \begin{equation}\label{eq:VB}
    \alpha(x) \equiv_{\FSL} \st_x(\phi_1) \cvee \cdots \cvee \st_x(\phi_n).
  \end{equation}
  If $n = 0$ then $\alpha(x) \equiv_{\FSL} \neg(x = x)$.
  Suppose $n \geq 1$.
  We claim that $\alpha(x) \equiv_{\FSL} \st_x(\phi_1 \vee \cdots \vee \phi_n)$.
  It is easy to see that
  $\st_x(\phi_1) \cvee \cdots \cvee \st_x(\phi_n)
    \models_{\FSL} \st_x(\phi_1 \vee \cdots \vee \phi_n)$, and hence
  $\alpha(x) \models_{\FSL} \st_x(\phi_1 \vee \cdots \vee \phi_n)$.
  So it suffices to prove
  $\st_x(\phi_1 \vee \cdots \vee \phi_n) \models_{\FSL} \alpha(x)$.
  
  Let $\fmo{M}$ be a $\FOL$-structure and $\fmo{M}_{\circ} = (W, 1, \fmeet, V)$
  and $w \in W$.
  Suppose $\fmo{M} \models \st_x(\phi_1 \vee \cdots \vee \phi_n)[w]$.
  Then $\fmo{M}_{\circ}, w \Vdash \phi_1 \vee \cdots \vee \phi_n$
  so by definition there exists $v_1, \ldots, v_n \in W$ such that
  $v_1 \fmeet \cdots \fmeet v_n \fleq w$ and
  $\fmo{M}_{\circ}, v_i \Vdash \phi_i$ for all $i \in \{ 1, \ldots, n \}$.
  Therefore $\fmo{M} \models \st_x(\phi_i)[v_i]$ for all $i \in \{ 1, \ldots, n \}$,
  so by~\eqref{eq:VB} we have
  $\fmo{M} \models \alpha(x)[v_i]$ for all $i \in \{ 1, \ldots, n \}$.
  Now Lemma~\ref{lem:alphax-filter} implies $\fmo{M} \models \alpha(x)[w]$.
\end{proof}

\section{Preventing preservation of classical contradictions}\label{sec:7b}

  The characterisation result in Theorem~\ref{thm:VB} requires manual treatment of the case
  where $\alpha(x)$ is equivalent to $\neg(x = x)$. This is necessary because the
  standard translation of $\bot \in \lan{L}$ is $\st_x(\bot) = \forall y(yRx)$,
  rather than the usual $\neg(x = x)$. This, in turn, is a consequence of
  having an inconsistent state in \Lone-models.
  In this final section we suggest a method to resolve this imperfection
  by adapting the definition of a meet-simulation.
  Rather than a relation between pairs of states from one model and single
  states of another  model, we relate finite subsets of one model to single
  states of the other.
  
  We denote by $\fun{P}_{\omega}W$ the collection of finite subsets of $W$,
  and call the adapted notion of a meet-simulation a
  meet$_{\omega}$-simulation.

\begin{definition}\label{def:meet-omega-sim}
  A \emph{meet$_{\omega}$-simulation}
  between $\mo{M} = (W, 1, \fmeet, V)$ and $\mo{M}' = (W', 1', \fmeet', V')$ is
  a relation $T$ between $\fun{P}_{\omega}W$ and $W'$,
  such that for all $(X, w') \in T$:
  \begin{enumerate}\itemindent=2em\itemsep=1pt
    \renewcommand{\labelenumi}{\textup{(\theenumi)}\enspace}
    \renewcommand{\theenumi}{M$'_{\arabic{enumi}}$}
    \item \label{it:msim-1p}
          If $w \in V(p)$ for all $w \in X$, then $w' \in V'(p)$, for each $p \in \Prop$;
    \item \label{it:msim-2p}
          If $w = 1$ for all $w \in X$, then $w' = 1'$;
    \item \label{it:msim-3p}
          If for each $w \in X$ the elements $u_w, v_w \in W$ are such that
          $u_w \fmeet v_w \fleq w$, then there exist $v', u' \in W'$ such that
                  $(\{ u_w \mid w \in X \}, u') \in T$ and
                  $(\{ v_w \mid w \in X \}, v') \in T$
                  and $v' \fmeet' u' \fleq' w'$.
  \end{enumerate}
\end{definition}

  Meet-simulations satisfy each of these conditions, so we have:

\begin{lemma}\label{lem:meet-to-omega}
  Every meet-simulation is a meet$_{\omega}$-simulation.
\end{lemma}

  A routine induction on the structure of $\phi$ proves the following
  analogue of Propositions~\ref{prop:adeq} and~\ref{prop:adeq-meet-sim}.

\begin{proposition}\label{prop:adeq-omega-sim}
  Let $T$ be a meet$_{\omega}$-simulation between $\mo{M} = (W, 1, \fmeet, V)$
  and $\mo{M}' = (W', 1', \fmeet', V')$.
  Then for each $(X, w') \in T$ and every formula $\phi \in \lan{L}$:
  $$
    \text{if}\quad \mo{M}, w \Vdash \phi \text{ for all $w \in X$, }
    \quad\text{then}\quad \mo{M}', w' \Vdash \phi.
  $$
\end{proposition}

\begin{definition}
  A $\FOL$-formula $\alpha(x)$ is said to be preserved by meet$_{\omega}$-simulations
  if for every meet$_{\omega}$-simulation $T$ between models $\mo{M}$ and $\mo{M}'$
  and all $(X, w') \in T$ we have:
  $$
    \text{if}\quad \mo{M} \models \alpha(x)[w] \text{ for all $w \in X$,}
    \quad\text{then}\quad \mo{M}' \models \alpha(x)[w'].
  $$
\end{definition}

  The crucial observation of this section is that classical
  contradictions are not preserved by meet$_{\omega}$-simulations.
  This is shown in the following example.

\begin{example}\label{exm:pres-falsum}
  Let $\beta(x) \in \FOL$ be a contradiction, such as $\neg(x = x)$.
  Consider any \Lone-model $\mo{M} = (W, 1, \fmeet, V)$.
  Then $T = \{ (\emptyset, 1) \}$ is a meet$'$-simulation on $\mo{M}$.
  We vacuously have $\mo{M} \models \beta(x)[w]$ for every $w \in \emptyset$,
  but not $\mo{M} \models \beta(x)[1]$.
  So $\beta(x)$ is not preserved by meet$_{\omega}$-simulations.
\end{example}

\begin{theorem}
  Let $\alpha(x) \in \FOL$ be a formula with one free variable $x$.
  Then $\alpha(x)$ is equivalent over $\cat{FSL}$ to
  the standard translation of an $\lan{L}$-formula iff it is
  preserved by meet$_{\omega}$-simulations.
\end{theorem}
\begin{proof}
  The direction from left to right follows from combining 
  Propositions~\ref{prop:st1} and~\ref{prop:adeq-omega-sim}.
  For the converse, suppose that $\alpha(x)$ is preserved by meet$_{\omega}$-simulations.
  Since every meet-simulation is a meet$_{\omega}$-simulation by
  Lemma~\ref{lem:meet-to-omega}, $\alpha(x)$ is also preserved by meet-simulations.
  It then follows from Theorem~\ref{thm:VB} that $\alpha(x)$ is
  either a contradiction or equivalent over $\cat{FSL}$ to the standard
  translation of a $\lan{L}$-formula.
  Since Example~\ref{exm:pres-falsum} shows that $\alpha(x)$ cannot be a contradiction,
  the theorem follows.
\end{proof}

\begin{remark}
  Inspection of the proofs above shows that we can replace
  $\fun{P}_{\omega}W$ in Definition~\ref{def:meet-omega-sim}
  by $\fun{P}_{\leq 2}W$, the collection of subsets of $W$ containing at most
  two states.
\end{remark}

\section{Conclusion}\label{sec:8}

  We have characterised non-distributive positive logic as a fragment of
  first-order logic with one binary predicate and unary a predicate for
  each proposition letter. We first used \Lone-simulations to approach the
  characterisation. Subsequently, we observed that classical disjunctions of
  standard translations of formulae are still preserved by \Lone-simulations.
  To remedy this we introduced meet-simulations between models, which relate
  pairs of elements of a model to single states of a second model.
  These allowed us to give a Van Benthem style characterisation theorem
  for non-distributive first-order logic.
  
  It would be worth investigating meet-simulations further.
  In particular, it is not clear how the composition of meet-simulations
  should be defined and whether this results in a meet-simulation again.
  We also wonder if the largest meet-simulation on a model yields
  a congruence, hence a quotient~model.
  
  Another potential avenue for further research is to extend the results in
  this paper to other logics.
  Obvious choices are modal extensions of non-distributive lattice logic,
  such as the one studied in~\cite{Dmi21,BezEA22}.
  Additionally, there may be connections with (bi)simulations for team semantics
  of dependence logics~\cite{HelEA14,KinEA14,YanVaa17}
  and modal information logic~\cite{Ben22}.

\paragraph{Acknowledgements.}
I am grateful to the anonymous reviewer for their insightful comments,
and to Anna Dmitieva for pointing out several typos and an error
in a draft of the manuscript.

\footnotesize{
\bibliographystyle{plainnat}
\bibliography{sim-lattice-logic.bib}
}

\end{document}